\documentclass[oneside,a4paper]{article}
\usepackage{lp}

\begin{document}
\title{Alexandrov meets Lott--Villani--Sturm}
\author{Anton Petrunin}
\date{}
\maketitle 

\begin{abstract}
Here I show compatibility of two definition of generalized curvature bounds ---
the lower bound for sectional curvature in the sense of Alexandrov and lower bound for Ricci curvature in the sense of Lott--Villani--Sturm.
\end{abstract}

\section*{Introduction}

Let me denote by $\CD m\kappa$ the class of  metric-measure spaces which satisfy a weak curvature-dimension condition for dimension $m$ and curvature $\kappa$ (see preliminaries).
By $\CBB m\kappa$, I will denote the class of all $m$-dimensional Alexandrov spaces with curvature $\ge\kappa$ equipped with the volume-measure (so $\CBB m\kappa$ is a class of metric-measure spaces).
\bigskip

\parbf{Main theorem.}
$\CBB m0\i \CD m{0}$.
\bigskip

The question appears first in \cite[7.48]{lott-villani}.
In \cite{villani}, it is formulated more generally: 
$\CBB m\kappa\i \CD m{(m-1)\kappa}$.
The later statement can be proved, 
along the same lines, 
but I do not write it down.

\parbf{About the proof.} 
The idea of the proof is the same as in the  Riemannian case (see \cite[6.2]{CMS} or \cite[7.3]{lott-villani}).
One only needs to extend certain calculus to Alexandrov spaces.
To do this, I use the technique developed in \cite{petrunin:HarmFun}.
Let me illustrate the idea on a simpler problem.

Let $M$ be a $2$-dimensional non-negatively curved Riemannian manifold and $\gamma_\tau\:[0,1]\to M$ be a continuous family of unit-speed geodesics such that
$$|\gamma_{\tau_0}(t_0)-\gamma_{\tau_1}(t_1)|\ge|t_1-t_0|.
\eqlbl{equi-dist}$$
for any $t_0,t_1,\tau_0,\tau_1$ (here $|x-y|$ denotes the distance between points $x$ and $y$ in a metric space). 
Set $\ell(t)$ to be the total length of curve $\sigma_t\:\tau\mapsto\gamma_\tau(t)$.
Then $\ell(t)$ is a concave function --- that is easy to prove.

Now, assume you have an Alexandrov space $A\in\CBB 20$ instead of Riemannian manifold $M$ and a non-continuous family of unit-speed geodesics $\gamma_\tau(t)$ which satisfies \ref{equi-dist}.
Define $\ell(t)$ as the 1-dimensional Hausdorff measure of image of $\sigma_t$.
Then  $\ell$ is also concave.

Here is an idea how one can proceed; 
it is not the simplest one but the one which admits a proper generalization.
Consider two distance functions $\psi=\dist_{\Im\sigma_0}$ and $\phi=\dist_{\Im\sigma_1}$.
Note that geodesics $\gamma_\tau(t)$ are also gradient curves of $\psi$ and $\phi$.
This implies that $\Lap\phi+\Lap\psi$ vanishes almost everywhere on the image of the map $(\tau,t)\to\gamma_\tau(t)$ (Laplacian of semiconcave function on Alexandrov space is a Radon sign-measures).
Then the result follows from the second variation formula from \cite{petrunin:parallel}
and calculus on Alexandrov spaces developed in \cite{perelman:DC}.

\parbf{Remark.}
Although $\CD m \kappa$ is a very natural class of metric-measure spaces, 
some basic tools in Ricci comparison can not work there in principle.
For instance,
\emph{there are $\CD m 0$-spaces which do not satisfy the Abresch--Gromoll inequality}, (see \cite{AG}).
Thus, one has to modify the definition of the class $\CD m \kappa$
to make it suitable for substantial applications in Riemannian geometry.

\smallskip

I'm grateful to 
V.~Kapovitch,
A.~Lytchak 
and C.~Villani, 
for their help.

\section{Preliminaries}

\parbf{Prerequisite.} 
The reader is expected to be familiar with 
basic definitions and notions of optimal transport  theory as in \cite{villani}, 
measure theory on Alexandrov spaces from \cite{BGP}, 
DC-structure on Alexandrov's spaces from \cite{perelman:DC} 
and technique and notations of gradient flow as in \cite{petrunin:survey}.

\parbf{What needs to be proved.} 
Let me recall the definition of class $\CD m0$ only --- it is sufficient for understanding this paper.
The definition of $\CD m\kappa$ can be found in \cite[29.8]{villani}.

Similar definitions were given in \cite{lott-villani} and \cite{sturm}.
The idea behind these definitions --- convexity of certain functionals in the Wasserstein space
over a Riemannian manifold, appears in \cite{otto-villani}, \cite{CMS}, \cite{SvR}.
In the Euclidean context, this notion of convexity goes back to \cite{McC}.
More on the history of the subject can be found in \cite{villani}.

For a metric-measure space $X$, 
I will denote 
by $|x-y|$ the distance between points $x,y\in X$
and by $\vol E$ the distinguished measure of Borel subset $E\i X$ (I will call it \emph{volume}).
Let us denote by $\Ptwo X$ the set of all probability measures with compact support in $X$ equipped with Wasserstein distance of order 2, see \cite[6.1]{villani}.

Further, we assume $X$ is a proper geodesic space;
in this case $\Ptwo X$ is geodesic.

Let $\mu$ be a probability measure on $X$.
Denote by $\mu^{r}$ the absolutely continuous part of $\mu$ with respect to volume.
That is, 
$\mu^{r}$ coincides with $\mu$ outside a Borel subset of volume zero and
 there is a Borel function $\rho\:X\to\RR$ such that $\mu^{r}\z=\rho{\cdot}\vol$.
Define
$$U_m\mu
\df
\int\limits_X\rho^{1-\frac1m}{\cdot}\vol
=
\int\limits_X\frac1{\sqrt[m]{\rho}}{\cdot}\mu^r.$$

Then $X\in \CD m0$ if the functional $U_m$ is concave on  $\Ptwo X$;
that is, for any two measures $\mu_0,\mu_1\in \Ptwo X$,
there is a geodesic path $\mu_t$, in $\Ptwo X$, $t\in[0,1]$ such that the real function $t\mapsto U_m\mu_t$ is concave.

\parbf{Calculus in Alexandrov spaces.} 
Let $A\in\CBB m\kappa$ and $S\i A$ be the subset of singular points;
that is, $x\in S$ iff its tangent space $\T_x$ is not isometric to Euclidean $m$-space $\EE^m$. 
The set $S$ has zero volume (\cite[10.6]{BGP}).
The set of regular points $A\backslash S$ is convex (\cite{petrunin:parallel}); 
that is, any geodesic connecting two regular points contains only regular points.

According to \cite{perelman:DC},  
if $f\:A\to\RR$ is a semiconcave function 
and $\Omega\i A$ is an image of a $\mathrm{DC}_0$-chart,
then
$\partial_k f$ and components of metric tensor $g^{i j}$ are functions of locally bounded variation which are continuous in $\Omega\backslash S$.

Further, for almost all $x\in A$ the Hessian of $f$ is defined.
That is, there is a subset of full measure $\Reg f\i A\backslash S$ such that for any $p\in \Reg f$ there is a bi-linear form%
\footnote{Since $p\in A\backslash S$, its tangent space $\T_p$ is isometric to Euclidean $m$-space. Therefore we can talk about bi-linear forms on $\T_p$.}
 $\Hess_p f$ on $\T_p$ such that 
$$f(q)=f(p)+d_p f(v)+\tfrac12{\cdot}[\Hess_p f](v,v)+o(|v|^2),$$
where $v=\log_p q$.
Moreover, the Hessian can be found using standard calculus in the $\mathrm{DC}_0$-chart.
In particular,
$$\Trace[\Hess f]
\ae
\frac{\partial_i(\det g\cdot g^{i j}\cdot\partial_j f)}{\det g}.$$

Let us reformulate the second variation formula \cite[1.1B]{petrunin:parallel}
using of ultrafilters.
Let $\o$ be a nonprinciple ultrafiler on natural numbers,
$A\in\CBB{m}{0}$
and $[p q]$ be a minimizing geodesic in $A$ which is extendable beyond $p$ and $q$.
Assume further that one of (and therefore each) of the  points $p$ and $q$ is regular.
Then there is a model configuration $\~p,\~q\in\EE^m$ and isometries $\imath_p\:\T_p A\to\T_{\~p}\EE^m$,
$\imath_q\:\T_q A\to\T_{\~q}\EE^m$ 
such that for any fixed $v\in\T_p$ and $w\in\T_q$ we have
$$
\l|\exp_p(\tfrac 1 n{\cdot} v)\,\exp_q(\tfrac 1 n{\cdot} w)\r|
\le
\l|\exp_{\~p}\circ\imath_p(\tfrac 1 n{\cdot} v)\,\exp_{\~q}\circ\imath_q(\tfrac1 n{\cdot} w)\r|+o(n^2)$$
for $\o$-almost all $n$ (once the left-hand  side is defined).

{\sloppy If $\~\tau\:\T_{\~p}\to\T_{\~q}$ is the parallel translation in $\EE^m$, then the isometry  \hbox{$\tau\:\T_p\to\T_q$} which satisfy identity $\imath_q\circ\tau=\~\tau\circ\imath_p$ will be called a  \emph{parallel transportation} from $p$ to $q$.

}
\parbf{Laplacians of semiconcave functions.}
Here are some facts from \cite{petrunin:HarmFun}.

Given a function $f\:A\to\RR$, define its \emph{Laplacian} $\Lap f$ to be a Radon sign-measure which satisfies the following identity
$$\int\limits_Au{\cdot}\Lap f
=
-\int\limits_A\<\nabla u,\nabla f\>{\cdot}\vol$$
for any Lipschitz function $u\:A\to\RR$.

\begin{thm}{Claim}\label{lap}
Let $A\in\CBB{m}{\kappa}$ and $f\:A\to \RR$ be $\lambda$-concave Lipschitz function.
Then Laplacian $\Lap f$ is defined and
$$\Lap f\le m{\cdot}\lambda{\cdot}\vol.$$
In particular, $\Lap^s f$ --- the singular part $\Lap f$ is nonpositive.

Moreover, 
$$\Lap f =\Trace[\Hess f]{\cdot}\vol+\Lap^s f.$$

\end{thm}

\parit{Proof.} 
Let us denote by $F_t\:A\to A$ the $f$-gradient flow for time $t$.

Given a Lipschitz function $u\:A\to\RR$, 
consider family $u_t(x)=u\circ F_t(x)$.
Clearly, $u_0\equiv u$ and $u_t$ is Lipschitz for any $t\ge0$.
Further, for any $x\in A$ we have $\bigl|{\tfrac{d^+}{d t}}u_t(x)|_{t=0}\bigr|\le\Const$.
Moreover 
$$\tfrac{d^+}{d t}u_t(x)|_{t=0}\ae d_x u(\nabla_x f)\ae\<\nabla_x u,\nabla_x f\>.$$

Further,
$$\int\limits_A u_t{\cdot}\vol=\int\limits_A u{\cdot}(F_t\#\vol),$$ 
where $\#$ stands for push-forward. 
Since $|F_t(x)-F_t(y)|\le e^{\lambda t}{\cdot}|x-y|$  
(see \cite[2.1.4(i)]{petrunin:survey}), 
for any $x,y\in A$ we have
$$F_t\#\vol\ge \exp(-m{\cdot}\lambda{\cdot} t){\cdot}\vol.$$
Therefore, for any non-negative Lipschitz function $u\:A\to\RR$,
$$\int\limits_A u_t{\cdot}\vol
=
\int\limits_A u{\cdot} (F_t\#\vol)
\ge
\exp(-m{\cdot}\lambda{\cdot} t)\cdot\int\limits_A u{\cdot}\vol.$$
Therefore
$$\int\limits_A\<\nabla u,\nabla f\>{\cdot}\vol 
=
\l.\tfrac{d^+}{d t}\int\limits_A u_t{\cdot}\vol\r|_{t=0}
\ge
-m{\cdot}\lambda{\cdot}\int\limits_A u{\cdot}\vol.$$
That is, there is a Radon measure $\chi$ on $A$, such that
$$\int\limits_A u{\cdot}\chi
=
\int\limits_A\l[\<\nabla u,\nabla f\>+m{\cdot}\lambda{\cdot} u\r]{\cdot}\vol.$$
Set $\Lap f=-\chi+m{\cdot}\lambda$, 
it is a Radon sign-measure and $\chi=-\Lap f+m{\cdot}\lambda\ge 0$.

To prove the second part of theorem,
assume $u$ is a non-negative Lipschitz function with support in a DC$_0$-chart $U\to A$,
where $U\i \RR^m$ is an open subset.
Then 
\begin{align*}
\int\limits_A\<\nabla u,\nabla f\>
&=
\int\limits_U \det g\cdot g^{i j}\cdot\partial_i u\cdot\partial_j f \cdot 
d x^1{\cdot}d x^2\cdots d x^m
=
\\
&=-\int\limits_U u\cdot\partial_i(\det g\cdot g^{i j}\cdot\partial_j f) \cdot d x^1{\cdot}d x^2\cdots d x^m,
\end{align*}
Thus 
$$\Lap f
=
\partial_i(\det g\cdot g^{i j}\cdot\partial_j f)\cdot d x^1{\cdot}d x^2\cdots d x^m
\ae \Trace[\Hess f].$$
\qedsf

\parbf{Gradient curves.} Here I extend the notion of gradient curves to  families of functions, see \cite{petrunin:survey} for all necessary definitions.

Let $\II$ be an open real interval
and $\lambda\:\II\to\RR$ be a continuous function.
A one parameter family of functions $f_t\:A\to \RR$, $t\in \II$ will be called \emph{$\lambda(t)$-concave} if the function $(t,x)\mapsto f_t(x)$ is locally Lipschitz 
and $f_t$ is $\lambda(t)$-concave for each $t\in \II$.

A locally Lipschitz curve $\alpha\:\II\to A$ will be called an $f_t$-gradient curve
if for any $t\in \II$ and minimizing geodesic $[\alpha(t)p]$ 
\[|p-\alpha(t+\eps)|\le |p-\alpha(t)|-d_{\alpha(t)}f_t(\dir {\alpha(t)}p)\cdot\eps+o(\eps),\]
for $\eps>0$.
Note that if the right derivative $\alpha^+(t)$ is defined, then the last inequality implies that
\[\langle\alpha^+(t),v\rangle\ge d_{\alpha(t)}f(v)\]
for any $v\in\T_p$.
In fact the latter property is equivalent to the definition;
this follows since $\alpha$ is differentiable almost everywhere.

For a time-independent family of functions $f_t=f$ this condition is equivalent to the standard definiton of gradient flow by $\alpha^+(t)=\nabla_{\alpha(t)} f$.
Note however that the last identity implies that $\langle\alpha^+(t),\alpha^+(t)\rangle=d_{\alpha(t)}(\alpha^+(t))$, but analogous identity for a family of functions dose not hold.
For example, the curve $\alpha(t)=t$ is a $f_t$-gradient curve for the family $f_t(x)=-|x-t|$ on the real line while $\nabla_t f_t=0\ne 1=\alpha'(t)$.

The following is a slight generalization of \cite[2.1.2$\&$2.2(2)]{petrunin:survey}; it can be proved along the same lines.

\begin{thm}{Proposition-definition}\label{prop-def}
Let $A\in\CBB m\kappa$,
$\II$ be an open real interval, 
$\lambda\:\II\to\RR$ be a continuous function and 
$f_t\:A\to \RR$, $t\in \II$ be $\lambda(t)$-concave family.

Then for any $x\in A$ and $t_0\in \II$ there exists an $f_t$-gradient curve $\alpha$ which is defined in a neighborhood of $t_0$ and such that $\alpha(t_0)=x$.

Moreover, if $\alpha,\beta\:\II\to A$ are $f_t$-gradient then for any $t_0,t_1\in\II$, $t_0\le t_1$,
$$|\alpha(t_1)-\beta(t_1)|\le L{\cdot}|\alpha(t_0)-\beta(t_0)|,$$
where $L=\exp\l(\int_{t_0}^{t_1}\lambda(t){\cdot}\,d t\r).$
\end{thm}

Note that the above proposition implies that the value $\alpha(t_0)$ of an $f_t$-gradient curve $\alpha(t)$ uniquely determines it for all $t\ge t_0$ in $\II$. 
Thus we can define \emph{$f_t$-gradient flow} --- a family of maps $F_{{t_0},{t_1}}\: A\to A$ such that
$$F_{{t_0},{t_1}}\:\alpha(t_0)\mapsto\alpha(t_1),$$
where $\alpha$ is an $f_t$-gradient curve.

\begin{thm}{Claim}\label{vol-lap}
{\sloppy Let $f_t\:A\to\RR$ be a $\lambda(t)$-concave family and $F_{{t_0},{t_1}}$ be $f_t$-gradient flow.
Let $E\i A$ be a bounded Borel set. 
Fix $t_1$ and consider the function 
\hbox{$v(t)=\vol F^{-1}_{{t},{t_1}}(E)$}.
Then 
$$\bigl.v\bigr|_t^{t_1}
=
\int\limits_t^{t_1}\Lap f_\under{t}\l[F^{-1}_{\under{t},{t_1}}(E)\r]{\cdot}\,d\under{t}.$$

}
\end{thm}

\parit{Proof.}
Let $u\:A\to\RR$ be a Lipschitz function with compact support.
Set $u_t=u\circ F_{{t},{t_1}}$.
Note that the function $(x,t)\mapsto u_t(x)$ is locally Lipschitz.
Thus, the function 
$$w_u\:t\mapsto \int\limits_A u_t{\cdot}\vol$$ 
is locally Lipschitz.
Further
$$w_u'(t)\ae
-\int\limits_A \<\nabla u_t,\nabla f_t\>{\cdot}\vol
=
\int\limits_A u_t{\cdot}\Lap f_t.$$
Therefore 
$$\bigl.w_u\bigr|_t^{t_1}
=
\int\limits_t^{t_1}d\under{t}\cdot\int\limits_A u_{\under{t}}{\cdot}\Lap f_{\under{t}}.$$
The last formula extends to an arbitrary Borel function $u\:A\to\RR$ with bounded support.
Applying it to the characteristic function of $E$ we get the result.
\qeds

\section{Games with Hamilton--Jacobi shifts.}

Let $A\in\CBB m 0$.
For a function $f\:A\to \RR\cup\{+\infty\}$, let us define its Hamilton--Jacobi shift\footnote{There is a lot of similarity between the Hamilton--Jacobi shift of a function and an equidistant for a hypersurface.} $\HJ_{t}f\:A\to\RR$ for time $t>0$ as follows
$$(\HJ_{t}f)(x)\df\inf_{y\in A} \l\{f(y)+\tfrac1{2\cdot t}{\cdot}|x-y|^2\r\}.$$
We say that $\HJ_{t}f$ is defined if the above infimum is $>-\infty$ everywhere in $A$.
Note that
$$\HJ_{t_0+t_1}f=\HJ_{t_1}\HJ_{t_0}f,
\eqlbl{ding202}$$
for any $t_0,t_1>0$.
(The inequality $\HJ_{t_0+t_1}f\le\HJ_{t_1}\HJ_{t_0}f$ is a direct consequence of triangle inequality and it is actually equality for any intrinsic metric, in particular for Alexandrov space.)

Note that for $t>0$, $f_t=\HJ_{t}f$ forms a $\tfrac1t$-concave family,
thus, we can apply \ref{prop-def} and \ref{vol-lap}.
The next claim gives a more delicate property of the gradient flow for such families;
it is an analog of \cite[3.3.6]{petrunin:survey}.

\begin{thm}{Claim}\label{clm:anti-lip}
Let $A\in\CBB m 0$ and
$f_0\:A\to\RR$ be function;
suppose $f_t=\HJ_t f_0$ is defined for $t\in(0,1)$.
Assume $\gamma\:[0,1]\to A$ is a geodesic path which is an $f_t$-gradient  curve for  $t\in(0,1)$
and $\alpha\:(0,1)\to A$ is another $f_t$-gradient curve.
If for some $t_0\in(0,1)$, $\alpha(t_0)=\gamma(t_0)$, then $\alpha(t)=\gamma(t)$ for all $t\in(0,1)$. 
\end{thm}

\begin{wrapfigure}{r}{45mm}
\begin{lpic}[t(-2mm),b(0mm),r(0mm),l(0mm)]{pics/alpha-gamma(0.40)}
\lbl[r]{3,23;$x$} 
\lbl[t]{47,20;$z$} 
\lbl[l]{108,23;$y$} 
\lbl[b]{54,45;$p$} 
\lbl[rb]{49,33;$\ell$} 
\lbl[tl]{87,49;$\alpha$} 
\lbl[t]{80,22;$\gamma$} \lbl[t]{72,72;$\sigma$}
\end{lpic}
\end{wrapfigure}

\parit{Proof.}
Set $\ell=\ell(t)=|\alpha(t)-\gamma(t)|$; it is locally Lipschitz function defined on $(0,1)$.
We have to show that if $\ell(t_0)=0$ for some $t_0$ then $\ell(t)=0$ for all $t$.

According to \ref{prop-def}, 
$\ell(t)=0$ for all $t\ge t_0$. 
In order to prove that $\ell(t)=0$ for all $t\le t_0$, it is sufficient to show that
$$\ell'\ge- [\tfrac 1 t+\tfrac{2}{1-t}]{\cdot}\ell$$ 
for almost all $t$.

Since $\alpha$ is locally Lipschitz, for almost all $t$, $\alpha^+(t)$ and $\alpha^-(t)$ are defined and \emph{opposite}%
\footnote{That is, $|\alpha^+(t)|=|\alpha^-(t)|$ and $\mangle(\alpha^+(t),\alpha^-(t))=\pi$}
to each other.

Fix such $t$ and set $x=\gamma(0)$, $z=\gamma(t)$, $y=\gamma(1)$, $p=\alpha(t)$, so $\ell(t)=|p-z|$.
Note that function
$$f_t+\tfrac1{2\cdot(1-t)}{\cdot}\dist^2_y
\eqlbl{ding203}$$
has a minimum at $z$. 
Extend a geodesic $[z p]$ by a both-sides infinite unit-speed quasigeodesic%
\footnote{A careful proof of existence of quasigeodesics can be found in \cite{petrunin:survey}.}
 $\sigma\:\RR\to A$,
so $\sigma(0)=z$ and $\sigma^+(0)=\dir{z}{p}$.
The function $f_t\circ\sigma\:\RR\to\RR$ is $\tfrac1t$-concave and from \ref{ding203},
$$f_t\circ\sigma(s)\ge f_t(z)+\<\gamma^+(t),\dir z p\>{\cdot} s-\tfrac1{2\cdot(1-t)}{\cdot} s^2.$$
It follows that 
\begin{align*}
\<\alpha^+(t),\sigma^+(\ell)\> 
&\ge 
d_pf_t(\sigma^+(\ell))=\\
&=
(f_t\circ\sigma)^+(\ell)\ge\\
&\ge
\<\gamma^+(t),\dir z p\>
- [\tfrac1t+\tfrac2{1-t}]{\cdot}\ell.
\end{align*}
Now,
\begin{enumerate}
\item Vectors $\sigma^\pm(\ell)$ are polar; therefore
$\<\alpha^\pm(t),\sigma^+(\ell)\>+\<\alpha^\pm(t),\sigma^-(\ell)\>\ge 0$.
\item Vectors $\alpha^\pm(t)$ are opposite; therefore $\<\alpha^+(t),\sigma^\pm(\ell)\>+\<\alpha^-(t),\sigma^\pm(\ell)\>= 0$.
\end{enumerate}
Since $\sigma^-(\ell)=\dir p z$, we get 
\[\<\alpha^+(t),\sigma^+(\ell)\>+\<\alpha^+(t),\dir p z\>= 0.\]
Therefore 
$$\ell'
=
-\<\alpha^+(t),\dir p z\>-\<\gamma^+(t),\dir z p\>
\ge
-[\tfrac1t+\tfrac2{1-t}]{\cdot}\ell.$$
\qedsf

\begin{thm}{Proposition}\label{HJ-2nd-var}
Let $A\in\CBB m0$,
$f\:A\to\RR$ be a  bounded continuous function and let $f_t=\HJ_t f$.
Assume $\gamma\:(0,a)\to A$ is an $f_t$-gradient curve which is also a constant-speed geodesic.
Assume that function
$$h(t)
\df
\Trace[\Hess_{\gamma(t)}f_t]$$
is defined for almost all $t\in(0,a)$. 
Then 
$$h'\le -\tfrac1m{\cdot} h^2$$
in the sense of distributions; 
that is, for any non-negative Lipschitz function $u\:(0,a)\z\to\RR$ with compact support
$$\int\limits_0^a \l(\tfrac1m{\cdot} h^2{\cdot}u-h{\cdot} u'\r){\cdot}\,d t\ge0.$$

\end{thm}

\parit{Proof.}
Since $h$ is defined a.e., all $\T_{\gamma(t)}$ for $t\in(0,a)$ are isometric to Euclidean $m$-space.
From \ref{ding202},
$$f_{t_1}(x)=\inf_{y\in A}\l\{f_{t_0}(y)+\frac{|x-y|^2}{2{\cdot}(t_1-t_0)}\r\}.$$
Thus, for a parallel transportation $\tau\:\T_{\gamma(t_0)}\to \T_{\gamma(t_1)}$ along $\gamma$,
we have 
$$[\Hess_{\gamma(t_1)}f_{t_1}](y,y)
\le
[\Hess_{\gamma(t_0)}f_{t_0}](x,x)
+\frac{|\tau (x)-y|^2}{t_1-t_0}$$
for any $x\in\T_{\gamma(t_0)}$ and $y\in\T_{\gamma(t_1)}$.
It remains to take the trace of the last formula with $y=c{\cdot}\tau(x)$ for appropriate $c>0$.
\qeds

\section{Proof of the main theorem}

Let $A\in\CBB m 0$; in particular $A$ is a proper geodesic space.
Let $\mu_t$ be a family of probability measures on $A$ for $t\in[0,1]$ which forms a \emph{geodesic path}\footnote{that is, constant-speed minimizing geodesic defined on $[0,1]$} in $\Ptwo A$ and both $\mu_0$ and $\mu_1$ are absolutely continuous with respect to volume on $A$.

Since Alexandrov's spaces are nonbranching,
\cite[30.32]{villani} implies that
it is sufficient to show that the function 
$$\Theta\:t\mapsto U_m\mu_t$$ 
is concave.

\def\Pii{\text{\scriptsize $\Pi$}}

According to \cite[7.22]{villani}, there is a probability measure $\Pii$ on the space of all geodesic paths in $A$ which satisfies the following condition:
If $\Gamma=\supp\Pii$ and $e_t\:\Gamma\to A$ is evaluation map  
$e_t\:\gamma\mapsto\gamma(t)$ then $\mu_t=e_t\#\Pii$.

The measure $\Pii$ is called the  \emph{dynamical optimal coupling} for $\mu_t$ and the measure $\pi=(e_0,e_1)\#\Pii$ is the corresponding \emph{optimal transference plan}.
The space $\Gamma$ will be considered further equipped with the metric
\[|\gamma-\gamma'|=\max_{t\in[0,1]}|\gamma(t)-\gamma'(t)|.\]

\parit{First we present $\mu_t$ as the push-forward  for gradient flows of two opposite families of functions.}
According to \cite[5.10]{villani}, there are optimal price functions
 $\phi,\psi\:A\to\RR$ such that
$$\phi(y)-\psi(x)\le\tfrac12{\cdot}|x-y|^2$$
for any $x,y\in A$ 
and equality holds for any $(x,y)\in\supp\pi$.
We can assume that $\psi(x)=+\infty$ for $x\not\in\supp\mu_0$ and $\phi(y)=-\infty$ for $y\not\in\supp\mu_1$.

Consider two families of functions
$$\psi_t=\HJ_t\psi\ \ \ \ \t{and}\ \ \ \ \phi_t=\HJ_{1-t}(-\phi).$$
Clearly, $\psi_t$ forms a $\tfrac1t$-concave family for $t\in (0,1]$ 
and $\phi_t$ forms%
\footnote{Note that usually $\phi_t$ is defined with opposite sign, but I wanted to work with semiconcave functions only.} 
a $\tfrac1{1-t}$-concave family for $t\in [0,1)$.

It is straightforward to check that for any $\gamma\in\Gamma$ and $t\in(0,1)$
$$\pm\<\gamma^\pm(t),v\>
=
d_{\gamma(t)}\psi_t(v)
=
-d_{\gamma(t)}\phi_t(v);$$
in particular,
$$\gamma^+(t)=\nabla\psi_t\ \ \ \ \t{and}\ \ \ \ \gamma^-(t)=\nabla\phi_t.
\eqlbl{ding204}$$

For $0<t_0\le t_1\le 1$, let us consider the maps $\Psi_{t_0,t_1}\:A\to A$ --- the gradient flow of $\psi_t$, defined by
$$\Psi_{t_0,t_1}\alpha(t_0)=\alpha(t_1),$$
where $t\mapsto \alpha(t)$ is a $\psi_t$-gradient curve.
Similarly, $0\le t_0\le t_1<1$, define map $\Phi_{t_1,t_0}\:A\to A$ 
$$\Phi_{t_1,t_0}\beta(t_1)=\beta(t_0),$$
where $1-t\mapsto \beta(t)$ is a $\phi_{1-t}$-gradient curve.
According to \ref{prop-def}, 
$$\Psi_{t_0,t_1}\ \t{is}\ 
\tfrac{t_1}{t_0}\t{-Lipschitz\ \ \ \ \  and}
\ \ \ \ \  
\Phi_{t_1,t_0}
\ \t{is}\ 
\tfrac{1-t_0}{1-t_1}\t{-Lipschitz.}
\eqlbl{ding205}$$

From \ref{ding204}, $e_{t_1}=\Psi_{t_0,t_1}\circ e_{t_0}$ and $e_{t_0}=\Phi_{t_1,t_0}\circ e_{t_1}$.
Thus, for any $t\in(0,1)$, the map $e_t\:\Gamma\to A$ is bi-Lipschitz.
In particular, for any measure $\chi$ on $A$, there is a uniquely determined one-parameter family of ``pull-back'' measures $\chi_t^*$ on $\Gamma$,
that is, 
such that $\chi^*_t E=\chi (e_t E)$ for any Borel subset $E\i \Gamma$.

Fix some $z_0\in(0,1)$ (one can take $z_0=\tfrac12$) and equip $\Gamma$ with the measure $\nu=\vol_{z_0}^*$.
Thus, from now on ``almost everywhere'' has sense in $\Gamma$, $\Gamma\times (0,1)$ and so on.

\parit{Now we will represent $\Theta$ in terms of families of functions on $\Gamma$.}
Note that $\mu_1=\Psi_{t,1}\#\mu_t$ and $\Psi_{t,1}$ is $\tfrac1t$-Lipschitz.
Since $\mu_1$ is absolutely continuous,
so is $\mu_t$ for all $t$.
Set  $\mu_t=\rho_t{\cdot}\vol$.
Note that from \ref{ding205}, we get that
$$\l(\frac{1-t_1}{1-t_0}\r)^m
\le
\frac{\rho_{t_1}(\gamma(t_1))}{\rho_{t_0}(\gamma(t_0))}
\le
\l(\frac{t_1}{t_0}\r)^m$$
for almost all $\gamma\in\Gamma$ and $0<t_0<t_1<1$.
For $\gamma\in\Gamma$ set $r_t(\gamma)=\rho_t(\gamma(t))$.
Then 
$$\Theta(t)
=
\int\limits_A \rho_t^{-\frac1m}{\cdot}\mu_t
=
\int\limits_\Gamma r_t^{-\frac1m}{\cdot}\Pii.
\eqlbl{ding206}$$
In particular, $\Theta$ is locally Lipschitz in $(0,1)$.

{\sloppy
\parit{Next we show that the  measure $\Lap\phi_t$ is absolutely continuous on $e_t\Gamma$ and that \hbox{$r_t(\gamma(t))=\rho_t(\gamma(t)){\cdot}\Lap\phi_t$} in a weak sense.}
From \ref{ding205}, $\vol_t^*=e^{w_t}{\cdot}\nu$ for some Borel function $w_t\:\Gamma\to\RR$.
Thus
$$\vol e_t E=\int\limits_{E}e^{w_t}{\cdot}\nu$$ 
for any Borel subset $E\i \Gamma$.
Moreover, for almost all $\gamma\in\Gamma$, we have that function $t\mapsto w_t(\gamma)$ is locally Lipschitz in $(0,1)$ (more precisely, $t\mapsto w_t(\gamma)$ coincides with a Lipschitz function outside of a set of zero measure).
In particular $\frac{\partial w_t}{\partial t}$ is defined a.e. in $\Gamma\times (0,1)$
and moreover 
$$w_t
\ae
\int\limits_{z_0}^{t}\frac{\partial w_\under{t}}{\partial \under{t}}{\cdot}\,d\under{t}.$$

}\medskip

Further, from \ref{clm:anti-lip}, if $0<t_0\le t_1<1$ then for any $\gamma\in\Gamma$,
$$\Psi_{t_0,t_1}(x)= \gamma(t_1)\ \ \Longleftrightarrow\ \ x=\gamma({t_0}),$$
$$\Phi_{t_1,t_0}(x)= \gamma({t_0})\ \ \Longleftrightarrow\ \ x=\gamma({t_1}).$$
Thus, for any Borel subset $E\i \Gamma$,
$$e_{t_1}E=\Psi_{t_0,t_1}\circ e_{t_0}E=\Phi^{-1}_{t_1,t_0}\l( e_{t_0}E\r),$$
$$e_{t_0}E=\Phi_{t_1,t_0}\circ  e_{t_1}E=\Psi^{-1}_{t_0,t_1}\l(  e_{t_1}E\r).$$
Set 
$$v(t)
\df
\vol e_t E
=
\int\limits_{E}e^{w_t}{\cdot}\nu.$$ 
From \ref{vol-lap},
$$v'(t)
\ae
\Lap\psi_t(e_t E)
\ae
-\Lap\phi_t(e_t E).$$ 
Thus,
$\Lap \psi_t+\Lap\phi_t=0$
everywhere on $e_t\Gamma$.
From \ref{lap},
$$\Lap\psi_t\le\tfrac m t{\cdot}\vol,\ \ \ \ \ \ \ \ \Lap\phi_t\le\tfrac m {1-t}{\cdot}\vol.$$
Thus, both restrictions  $\Lap \psi_t|_{e_t\Gamma}$ and $\Lap\phi_t|_{e_t\Gamma}$
are absolutely continuous with respect to volume.
Therefore
$$v'(t)
\ae
\int\limits_{e_t E}\Trace[\Hess\phi_t]{\cdot}\vol.
$$

For the one parameter family of functions $h_t(\gamma)=\Trace[\Hess_{\gamma(t)}\phi_t]$, we have
$$\bigl.v\bigr|_{z_0}^t
=
\int\limits_E(e^{w_t}-1){\cdot}\nu
=
\int\limits_{z_0}^t d\under{t}\cdot\int\limits_E h_\under{t}{\cdot}e^{w_\under{t}}{\cdot}\nu$$
or any Borel set $E\i \Gamma$. 
Equivalently, 
$$\frac{\partial w_t}{\partial t}
\ae 
h_t.$$
From \ref{HJ-2nd-var}, 
$$\frac{\partial h_t}{\partial t}\le -\tfrac1m{\cdot} h_t^2.$$
Thus, for almost all $\gamma\in\Gamma$, the following inequality holds in the sense of distributions:
$$\frac{\partial^2}{\partial t^2}\exp\l(\frac{w_t(\gamma)}m\r)
=
\l(\tfrac1{m^2}{\cdot}{h_t^2}+\tfrac1{m}{\cdot}\frac{\partial h_t}{\partial t}\r){\cdot}\exp\l(\frac{w_t(\gamma)}m\r)\le 0;
$$ 
that is, $t\mapsto \exp\l(\frac{w_t(\gamma)}m\r)$ is concave
--- more precisely, $t\mapsto \exp\l(\frac{w_t(\gamma)}m\r)$ coincides with a concave function almost everywhere.

Clearly, for any $t$ we have
$\mu=r_t{\cdot} e^{w_t}{\cdot}\nu$.
Thus, for almost all $\gamma$ there is a non-negative Borel function $a\:\Gamma\to\RR_\ge$ such that $r_t\ae a {\cdot}e^{-w_t}$.
Continue \ref{ding206},
$$\Theta(t)
=
\int\limits_\Gamma r_t^{-\frac1m}{\cdot}\Pii
=
\int\limits_\Gamma{e^{\frac{w_t}m}}{\cdot}{\sqrt[m]{a}}{\cdot}\Pii.
$$
That is, $\Theta$ is concave as an average of concave functions.
Again, more precisely, $\Theta$ coincides with a concave function a.e., but since $\Theta$ is locally Lipschitz in $(0,1)$ we get that $\Theta$ is concave. \qeds

\end{document}